\documentclass[final,leqno,onefignum,onetabnum]{article}

\usepackage{amsmath,amssymb}

\usepackage{dcolumn}
\usepackage{bm}

\usepackage[T1]{fontenc}
\usepackage[latin9]{inputenc}
\usepackage{amsmath,amssymb,subeqnar,overpic}
\usepackage{graphicx}
\usepackage{color}

\usepackage[paperheight=11in,paperwidth=8.5in]{geometry}

\newcommand*\samethanks[1][\value{footnote}]{\footnotemark[#1]}

\title{Classification of Spatio-Temporal Data via Asynchronous Sparse Sampling:  Application to Flow Around a Cylinder}

\author{Ido Bright\thanks{Department of Applied Mathematics, University of Washington, Seattle, WA. 98195-3925.  $^\ddagger$ ({kutz@uw.edu}). Questions, comments, or corrections to this document may be directed to that email address.} \and Guang Lin\thanks{Department of Mathematics and School of Mechanical Engineering, Purdue University, West Lafayette, IN 47907 } \and J. Nathan Kutz$^\ddagger$\samethanks[1]}

\begin{document}
\maketitle

\begin{abstract}
We present a novel method for the classification and reconstruction
of time dependent, high-dimensional data using sparse measurements,
and apply it to the flow around a cylinder. Assuming the data lies
near a low dimensional manifold (low-rank dynamics) in space and has periodic time dependency
with a sparse number of Fourier modes, we employ compressive sensing for accurately classifying the
dynamical regime.   We further show that we can reconstruct the full spatio-temporal behavior
with these limited measurements, extending previous results of compressive sensing that
apply for only a single snapshot of data.  The method can be used for building 
improved reduced-order models and designing sampling/measurement strategies that
leverage time asynchrony.
\end{abstract}



\section{Introduction}

Data-driven methods, often rooted in new innovations of machine learning~\cite{duda,bishop,murphy}, 
are becoming transformative tools in the study of complex dynamical systems~\cite{kutz:2013}.  
Such methods aim to take advantage of the ubiquitous observation across the physical, engineering
and biological sciences that meaningful input/output of macroscopic variables are encoded in low-dimensional patterns of dynamic 
activity~\cite{Cross:1993,jones07,rabinovich01,rabinovich08,laurent01,gold99,lowD2}.
Reduced-order models (ROMs) capitalize on this observation with the goal of producing real-time
computational schemes for modeling complex systems
(See the recent edited monograph addressing the state of the art in 
reduced order methods for computational modeling of complex parametrized systems~\cite{roms}).
Complimenting these efforts are new theoretical efforts using sparsity promoting techniques
\cite{candes2006robust,candes2006stable,candes2006near,donoho2006compressed,gilbert,baraniuk,baraniuk2,cssignal} which exploit this fact in order to deal with the practical constraints of a limited number of sensors
for data acquisition in the complex system~\cite{bright13,brunton13,gappy,Ganguli:2012,Field:2004}.  
Such findings lead us to conjecture that a principled approach
to sensing in such complex systems will allow us to probe the latent,
low-dimensional structure of dynamic activity relating to the fundamental
quantities of interest.  A direct implication of such a data-driven
modeling strategy is that we can gain traction on prediction, state-estimation
and control of the complex system. To this end, we demonstrate that
asynchronous sparse measurements in both time and space allows for a robust
method, i.e. a spatio-temporal compressive sensing method, for classifying
the dynamical behavior of a given system.   Once classified, it can be used for 
full state reconstruction and/or future state prediction of the full spatio-temporal
system, extending previous works~\cite{bright13,brunton13} which only reconstructed
single time snapshots of the data.

To demonstrate out methodology, we consider a practical implementation
of our technique on time-varying fluid flows that are ubiquitous in
modern engineering and in the life sciences. Specifically, we use
the asynchronous sparse measurements in time and space on the pressure
field around a cylinder to classify the Reynolds number and reconstruct the
full pressure field. This problem
was previously and successfully considered~\cite{bright13,brunton13} with
the view of only considering spatial compressive sensing. The addition
of asynchronous time measurements is the new innovation considered here allowing for
a more robust architecture for performing the classification task.

Our method is motivated by bio-inspired engineering design principals~\cite{1,2,3,4}
whereby birds, bats, insects, and fish routinely harness limited measurements
of unsteady fluid phenomena to improve their propulsive efficiency,
maximize thrust and lift, and increase maneuverability \cite{5,8,6,10,7,9}.
Indeed, despite their limited number of computational resources and
spatially distributed, noisy sensors, they achieve robust flight dynamics
and control despite potentially large perturbations in flight conditions
(e.g. wind gusts, bodily harm). Such observations suggest the existence
of low-dimensional structures in the flow-field that are sparsely
sampled and which are exploited for robust control purposes.

The primary innovation here revolves around the compressive sensing ideology. Compressed
sensing allows the reconstruction of a signal using a small number
of measurements, based on the fact that the signal has a sparse representation
in an appropriate basis~\cite{candes2006robust,candes2006stable,candes2006near,donoho2006compressed,gilbert,baraniuk,baraniuk2}.
This allows the use of less measurements than the Nyquist\textendash{}Shannon
criterion requires~\cite{nyquist1928certain,citeulike:1584479}.
Although typical applications are found in areas spanning signal
processing to computer vision, our objective is to use the methodology
for the numerical solution of partial differential equations (PDE) and/or direct observations
of complex systems.  We present a specific example - the incompressible Navier-Stokes flow
around a cylinder. However, the framework presented can be applied
to a much broader set of equations. At the method's core is the fact that
most complex systems, including the flow around a cylinder, exhibit
an underlying low-dimensional structure in their dynamics, such an assumption
is at the heart of the ROM methodology.
Unlike the standard application of compressive sensing for image processing, for example,
where the sparse (wavelet) basis is already known, the optimal sparse
basis considered here is computed directly from PDE simulations (e.g. through a proper orthogonal decomposition, for instance). Previously,
the method was developed to classify the Reynolds number for flow around a
cylinder given a single snapshot measurement~\cite{bright13}. However, by using asynchronous
temporal measurements as well, the method can be greatly improved. Thus a
spatio-temporal version of the compressive sensing architecture is developed and advocated
here as the primary result.

The paper is outlined as follows:  In Sec.~2 we outline the underlying methods used for
the spatio-temporal dimensionality reductions, highlighting how low-rank structures in
both space and time are equally exploited.  Section~3 outlines the method used
for reconstruction of the solution given a sparse set of measurements.  Classification
of the dynamical regime of the system, which is based upon libraries of dynamical
modes, is demonstrated in Sec.~4.  Sections 2-4 in combination form the core innovation
of this manuscript, leading to the example problem of flow around a cylinder of Sec.~5.
A brief conclusion and outlook are given in Sec.~6.

\section{Dimensionality Reduction for Data Snapshots\label{sec:Dimension-Reduction-Method}}

Dimensionality reduction has a long history in regard to theoretical methods applied to complex
dynamical systems.  Before the advent of computers, dimensionality reduction often was achieved
through asymptotic reductions of the governing equations.  Indeed, many well characterized
theoretical examples of fluid dynamics are due to perturbative reductions of the full system in
some parameter regime, e.g. the low Reynolds number limit.  Such reductions allow for a
simplified model where known analytic methods can be applied to great effect.  These asymptotic
techniques also serve as the backbone for normal form reductions for characterizing the low-dimensional
embedding subspaces where bifurcations 
in complex systems occur.   With the 
advent of the computer and modern high-performance computing, these same kind of reductions
can be achieved through principled approaches such as the proper orthogonal decomposition~\cite{PODHolmes},
which is also known with minor variation as Principal Component Analysis (PCA)~\cite{Pearson:1901}, the Karhunen--Lo\`eve (KL)
decomposition, the Hotelling transform~\cite{hotellingJEdPsy33_1,hotellingJEdPsy33_2} and/or 
Empirical Orthogonal Functions (EOFs)~\cite{lorenzMITTR56}.  These methods identify the
low-dimensional subspaces and modal basis that optimally represent, in an $\ell_2$ sense, the dynamics of the system, thus forming the basis of the ROM architecture~\cite{roms}.

For our purposes, consider data that comes from a dynamical system of a single spatial variable $x$ of the form 
\begin{equation}
u_{t}=N\left(\nu,x,t,u,u_{x},u_{xx},\dots\right),
\label{eq:Nu}
\end{equation}
where $N(\cdot)$ determines the evolution of the complex system which can be a function
of space ($x\in[-L,L]$), time ($t\in[0,T]$) and linear and nonlinear terms involving $u(x,t)$ and its derivatives.
Additionally, the parameter $\nu$ is a bifurcation parameter that determines the dynamical
regime of the system.  Such systems are typically solved using a numerical discretization scheme
where the $k$th time point ($k=1,2,\cdots,m$) is given by
\begin{equation}
  u(x,t_k)\rightarrow {\bf u}_k=[u(x_1,t_k) \,\, u(x_2,t_k) \,\, \cdots \,\, u(x_n,t_k)]^T
  \label{eq:discreteu}
\end{equation}
with the discretization grid $\Delta x=x_{j+1}\!-\!x_j \ll 1$, $x_1=-L$, $x_n=L$ and $n\gg 1$.
Thus in a simulation, one has access
to all the data in a discretized form in space and time.  Moreover the data is high-dimensional
since a fine grid (of size $n$) is typically required for accurate solutions of the  partial differential equation.

For various values of the bifurcation parameter $\nu$ we assume we have access to the data, either
from computation or from dense measurements in time and space. We can arrange the full data
into a matrix whose columns are the full discretized state of the system, or snapshots, sampled evenly in time
\begin{equation}
   {\bf A}_D^{\nu} = \left[ {\bf u}_1 \,\, {\bf u}_2 \,\, \cdots \,\, {\bf u}_m    \right]
\end{equation}
where ${\bf A}_D^{\nu}\in \mathbb{R}^{n\times m}$.
In the application we wish to consider in this manuscript, we assume that
we have access only to a subset of the data, based on the position in the spatial grid $x_j$, i.e. we
have a limited number of spatial sensors available that are fixed to a specific spatial location. 
To be more precise, consider the
the projection matrix ${\bf P}_k$ whose $k$ nonzero entries determine the spatial locations for data 
sensors (let $k=3$ for illustrative purposes):
\begin{equation}
\label{eq:P2}
  {\bf P}_3 = \left[  \begin{array}{cccccccccc}  
     1 & 0          &  \cdots     &     &      & & & & \cdots  & 0 \\
     0 & \cdots   & 0             &  1 &  0  & \cdots & & & \cdots  & 0 \\
     0 & \cdots   &                &     & \cdots     & 0 & 1 & 0 &\cdots & 0
  \end{array}    \right]
\end{equation}
This example maps out the full state of the system onto 3 measurement locations 
\begin{equation}
\tilde{\bf u}= {\bf P}_3{\bf u}.
\end{equation}
More generally for $k$ measurement locations, we can construct the data matrix
\begin{equation}
   {\bf A}_S^{\nu} = \left[ \tilde{\bf u}_1 \,\, \tilde{\bf u}_2 \,\, \cdots \,\, \tilde{\bf u}_m    \right]
\end{equation}
where ${\bf A}_S^{\nu}\in \mathbb{R}^{k\times m}$.
Thus ${\bf A}_{S}^{\nu}$ is a subset of data generated from a small number sensors $k$ where $k\ll n$.
Our hypothesis is that the data satisfies the following conditions: 

\begin{description}

\item [{H1}] The data in $A_{S}^{\nu}$ has a low dimensional representation
in the phase space, namely, it can be approximated by a low dimensional
vector space while maintaining most of the energy (variance) in the $\ell_2$
sense. 
\end{description}

To extract the low-rank representation of the data, we apply a singular value decomposition (SVD) to the data matrix to obtain 
\begin{equation}
{\bf A}_{S}^{\nu}={\bf V}_{S}^{\nu}{\bf \Sigma}_{S}^{\nu}{\bf U}_{S}^{\nu} .
\label{eq:svd}
\end{equation}
The SVD decomposition can be used for producing a principled low-dimensional 
representation of the data.  Truncation is typically achieved by inspection of the diagonal singular 
value matrix ${\bf \Sigma}_{S}^{\nu}$.   Specifically, $d^{\nu}$ modes are selected that, for instance, 
represent 99\% of the total variance in the data.  For data with noise, recent work by
Gavish and Donoho~\cite{gavish} provides a principled truncation algorithm to account for 
the effects of noise.

In our analysis, we reduce the data to $d^{\nu}$-dimensions for a variety of dynamical regimes 
$\nu$.  This allows us to obtain the approximation to the data matrix 
\begin{equation}
\bar{\bf A}_{S}^{\nu}=\bar{\bf V}_{S}^{\nu}\bar{\bf \Sigma}_{S}^{\nu}\bar{\bf U}_{S}^{\nu},
\end{equation}
where the barred matrices are $d^{\nu}$-rank approximations to their unbarred counterparts, ultimately providing
the approximation to the full data
\begin{equation}
{\bf A}_{S}^{\nu}\approx\bar{\bf A}_{S}^{\nu}=\sum_{i=1}^{d^{\nu}}\sigma_{i}\mathbf{v}_{i}^{\nu}\mathbf{u}_{i}^{\nu}.
\label{eq:lowrank}
\end{equation}
where $\bar{\bf V}_{S}^{\nu}=\left[\mathbf{v}_{1}^{\nu},\dots,\mathbf{v}_{d^{\nu}}^{\nu}\right]$,
$\bar{\bf U}_{S}^{\nu}=\left[ \mathbf{u}_{1}^{\nu} \,\, \mathbf{u}_{2}^{\nu}  \,\, \dots \,\, \mathbf{u}_{d^{\nu}}^{\nu}\right]$,
and $\sigma_{i}$ are the corresponding diagonal elements $\bar{\bf \Sigma}_{S}^{\nu}$.

In addition to our assumption that (\ref{eq:Nu}) exhibits low-dimensional dynamics, there is a large
subclass of dynamical systems that also exhibit time dynamics that are nearly periodic.  More precisely,
in the application considered here, the time dynamics can be approximated by a small number of Fourier
modes, thus allowing for a second sparsification step.  This leads to the second hypothesis of the current work:

\begin{description}
\item [{H2}] Each row of the matrix $\bar{\bf U}_{S}^{\nu}$ can be approximated,
in the $\ell_{2}$ sense, by a small number of Fourier modes in time, namely
by cosines and sines. 
\end{description}

Given this second hypothesis, our goal now is to further approximate the full data matrix ${\bf A}_{S}^{\nu}$ by
a sparse number of Fourier modes.  Thus in the modal decomposition (\ref{eq:svd}), the matrix
containing the time dynamic modes ${\bf u}_k^\nu$ is approximated by
\begin{equation}
\mathbf{u}_{i}^{\nu}=\sum_{j=1}^{n_{i}^{\nu}} \left( a_{i,j}^{c,\nu}\mathbf{c}_{i,j}^{\nu}+a_{i,j}^{s,\nu}\mathbf{s}_{i,j}^{\nu} \right) ,
\label{eq:timesparse}
\end{equation}
where $n_{i}^{\nu}$ is the small (sparse) number of Fourier modes retained with
\begin{subeqnarray}
 &&  \mathbf{c}_{i,j}^{\nu}=\left[\cos\left(\omega_{i,j}^{\nu}t_{1}\right),\cos\left(\omega_{i,j}^{\nu}t_{2}\right),\dots,\cos\left(\omega_{i,j}^{\nu}t_{N}\right)\right]  \\
 &&  \mathbf{s}_{i,j}^{\nu}=\left[\sin\left(\omega_{i,j}^{\nu}t_{1}\right),\sin\left(\omega_{i,j}^{\nu}t_{2}\right),\dots,\sin\left(\omega_{i,j}^{\nu}t_{N}\right)\right] \, .
\end{subeqnarray}
and where $\omega_{i,j}^{\nu}$ is the angular frequency obtained
by the Fourier transform.

The low-rank approximation by the SVD (\ref{eq:lowrank}) and the sparse Fourier mode approximation (\ref{eq:timesparse})
can be combined to give the time-space reduction
\begin{equation}
{\bf A}_{S}^{\nu}=\sum_{i=1}^{d^{\nu}}\sigma_{i}\mathbf{v}_{i}\sum_{j=1}^{n_{i}^{\nu}}\left(a_{i,j}^{c,\nu}\mathbf{c}_{i,j}^{\nu}+a_{i,j}^{s,\nu}\mathbf{s}_{i,j}^{\nu}\right)=\sum_{i=1}^{d^{\nu}}\sum_{j=1}^{n_{i}^{\nu}}\left(\alpha_{i,j}^{c,\nu}\mathbf{v}_{i}\mathbf{c}_{i,j}^{\nu}+\alpha_{i,j}^{s,\nu}\mathbf{v}_{i}\mathbf{s}_{i,j}^{\nu}\right),
\end{equation}
where $\alpha_{i,j}^{c,\nu}=\sigma_{i}a_{i,j}^{c,\nu}$ and $\alpha_{i,j}^{s,\nu}=\sigma_{i}a_{i,j}^{s,\nu}$.

We now perform the final dimensionality-reduction step in our analysis.  Specifically, 
we define $\left(\alpha_{i,j}^{\nu}\right)^{2}=\left(\alpha_{i,j}^{s,\nu}\right)^{2}+\left(\alpha_{i,j}^{c,\nu}\right)^{2}$ and
choose $m^{\nu}$ of the largest $\alpha_{i,j}^{\nu}$ so that
their sum of squares approximates the energy up to a given threshold for truncation.  This allows us
to obtain the final low-rank approximation of the data matrix
\begin{equation}
{\bf A}_{S}^{\nu}\approx\sum_{k=1}^{m^{\nu}}\left(\alpha_{k}^{c,\nu}\mathbf{v}_{i_{k}}^{\nu}\mathbf{c}_{i_{k},j_{k}}^{\nu}+\alpha_{k}^{s,\nu}\mathbf{v}_{i_{k}}^{\nu}\mathbf{s}_{i_{k},j_{k}}^{\nu}\right).\label{eq:sparse_app}  
\end{equation}
Note that we do not need to store the complete vectors $\mathbf{c}_{i_{k},j_{k}}^{\nu}$
and $\mathbf{s}_{i_{k},j_{k}}^{\nu}$, but only the vectors $\mathbf{v}_{i_{k}}^{\nu}$,
and the frequencies $\mathbf{\omega}_{i_{k},j_{k}}^{\nu}$ obtained
by the Fast Fourier Transform.  Approximation (\ref{eq:sparse_app}) allows us to use a small number of
measurements in time and space to reconstruct an approximation to the entire data matrix.  
This is the primary innovation of this work.

\section{Reconstruction from Sparse Measurements}

The primary objective in this work is to use sparse measurements in
space and time to reconstruct the full state  ${\bf u}$ of (\ref{eq:Nu}) and (\ref{eq:discreteu}).
To this end, consider that we are given $m^{\nu}$ measurements 
\begin{equation}
\mathbf{p}=\left[{\bf p}_{1} \,\, {\bf p}_2 \,\, \cdots \,\, {\bf p}_{m^\nu}\right]
\label{eq:pi}
\end{equation}
at the space-time points
\begin{equation}
\left(t_{i},x_{i}\right)=\left(\tau_{i},\chi_{i}\right) \,\,\,\,\, \mbox{for} \,\,\,\,\, i=1,2,\cdots,m^{\nu}.
\end{equation}
Recall that we assumed that the measurements can be respresented
sparsely by the parameter $\nu$ in (\ref{eq:sparse_app}),
namely: 
\begin{equation}
{\bf p}_{i}=\sum_{k=1}^{m^{\nu}}\left[\mathbf{v}_{i_{k}}^{\nu}\right]_{\chi_{i}}\left(\cos\left(\omega_{i_{k},j_{k}}^{\nu}\left(\tau_{i}+\phi\right)\right)+\sin\left(\omega_{i_{k},j_{k}}^{\nu}\left(\tau_{i}+\phi\right)\right)\right),
\end{equation}
where $\left[{\bf v}\right]_{l}$ denotes the $l$'th element of the vector,
and $\phi$ is an unknown phase.

We relax this representation by folding in the unknown phase into additional amplitude parameters.  The new representation
takes the form
\begin{equation}
{\bf p}_{i}=\sum_{k=1}^{m^{\nu}}\left[\mathbf{v}_{i_{k}}^{\nu}\right]_{\chi_{i}}\left(A_{i_{k},j_{k}}^{\nu}\cos\left(\omega_{i_{k},j_{k}}^{\nu}\tau_{i}\right)+B_{i_{k},j_{k}}^{\nu}\sin\left(\omega_{i_{k},j_{k}}^{\nu}\tau_{i}\right)\right).\label{eq:reconstr}
\end{equation}
Note that for each $\omega_{i_{k},j_{k}}^{\nu}$, there are distinct
$A_{i_{k},j_{k}}^{\nu}$ and $B_{i_{k},j_{k}}^{\nu}$ that do not
depend on $i$.   Thus there are a total of $2m^{\nu}$ amplitude parameters that need to be determined.

The representation (\ref{eq:reconstr}) is for a single measurement.  The total number of measurements $m^{\nu}$, as
given by (\ref{eq:pi}), can be represented in matrix form by
 \begin{equation}
 \mathbf{p}={\bf \Phi}^{\nu}\mathbf{b},
 \end{equation}
where the odd and even columns of ${\bf \Phi}^{\nu}$ are  
$\left[{\bf \Phi}^{\nu}\right]_{i,2j-1}\!=\!\left[\mathbf{v}_{i_{k}}^{\nu}\right]_{\chi_{i}} \!\! \cos (\omega_{i_{j},j_{j}}^{\nu}\tau_{i} )$  and $\left[{\bf \Phi}^{\nu}\right]_{i,2j}\!=\!
\left[v_{i_{j}}\right]_{\chi_{i}} \!\! \sin (\omega_{i_{j},j_{j}}^{\nu}\tau_{i} )$ so that
\begin{equation}
{\bf \Phi}^{\nu}\!\!=\!\!\! \left[ \!\!\! \begin{array}{cccc}
\left[\mathbf{v}_{i_{1}}\right]_{\chi_{1}}\!\cos\!\left(\omega_{i_{1},j_{1}}^{\nu}\tau_{1}\!\right) & \!\left[\mathbf{v}_{i_{1}}\right]_{\chi_{1}}\!\sin\!\left(\omega_{i_{1},j_{1}}^{\nu}\tau_{1}\!\right) & \!\cdots\! & \!\left[\mathbf{v}_{i_{m^{\nu}}}\right]_{\chi_{1}}\!\sin\!\left(\omega_{i_{m^{\nu}},j_{m^{\nu}}}^{\nu}\tau_{1}\!\right)\\
\left[\mathbf{v}_{i_{1}}\right]_{\chi_{2}}\!\cos\!\left(\omega_{i_{1},j_{1}}^{\nu}\tau_{2}\!\right) & \!\left[\mathbf{v}_{i_{1}}\right]_{\chi_{2}}\!\sin\!\left(\omega_{i_{1},j_{1}}^{\nu}\tau_{2}\!\right) & \!\cdots\! & \!\left[\mathbf{v}_{i_{m^{\nu}}}\right]_{\chi_{2}}\!\sin\!\left(\omega_{i_{m^{\nu}},j_{m^{\nu}}}^{\nu}\tau_{2}\!\right)\\
\vdots & \vdots & \!\ddots\! & \vdots\\
\left[\mathbf{v}_{i_{1}}\right]_{\chi_{N}}\!\cos\!\left(\omega_{i_{1},j_{1}}^{\nu}\tau_{N}\!\right) & \!\left[\mathbf{v}_{i_{1}}\right]_{\chi_{N}}\!\sin\!\left(\omega_{i_{1},j_{1}}^{\nu}\tau_{N}\!\right) & \!\cdots\! & \!\left[\mathbf{v}_{i_{m^{\nu}}}\right]_{\chi_{N}}\!\sin\!\left(\omega_{i_{m^{\nu}},j_{m^{\nu}}}^{\nu}\tau_{N}\!\right)
\end{array} \!\!\! \right]
\end{equation}
with 
\begin{equation}
{\bf b}=
\left[\begin{array}{c}
A_{i_{1},j_{1}}^{\nu}\\
B_{i_{1},j_{1}}^{\nu}\\
A_{i_{2},j_{2}}^{\nu}\\
B_{i_{2},j_{2}}^{\nu}\\
\vdots\\
B_{i_{m^{\nu}},j_{m^{\nu}}}^{\nu}
\end{array}\right].
\end{equation}
Note that the matrix ${\bf \Phi}^{\nu}$
depends on the time-space locations where measurements are made.  The reconstruction based
upon (\ref{eq:sparse_app}) can then be accomplished using a standard Moore-Penrose pseudo-inverse~\cite{trefethen} 
given the measurements ${\bf p}$.  This allows for computing the unknowns
${\bf b}$ and performing the full state reconstruction with limited measurements.

\section{Dynamical Libraries and Classification}

The previous section considered our ability to perform a reconstruction or approximation of the full state of the system
using a limited number of measurements.  The reconstruction was for a specific value of the bifurcation parameter $\nu$.
However, there may be many different dynamical regimes of (\ref{eq:Nu}) that are of interest and have low-dimensional
dynamics~\cite{bright13,brunton13}.   For each dynamical regime $\nu_j$, we can then construct the 
measurement matrices ${\bf \Phi}^{\nu_{1}},\dots, {\bf \Phi}^{\nu_{M}}$ where we assume there are $M$
dynamical regimes of interest ($\nu_1, \nu_2, \cdots , \nu_M$).  We can collect these various measurement matrices into a single library matrix:
\begin{equation}
{\bf \Psi}=\left[\Phi^{\nu_{1}} \,\, \Phi^{\nu_{2}} \,\, \dots \,\, \Phi^{\nu_{M}}\right].
\end{equation}
Thus ${\bf \Psi}$ contains all the low-rank approximations for space-time measurements from the $M$
dynamical regimes.
Such a database building strategy is also at the forefront of ROM architectures since
ROMs are neither robust with respect to parameter changes nor cheap to generate. Thus methods based on a database of ROMs coupled with a suitable interpolation schemes greatly reduces the computational cost for aeroelastic predictions while retaining good accuracy~\cite{roms,amsallem}.

Our objective ultimately is to use a small number of measurements to reconstruct the full state of
the system.  However, in order to do so efficiently and effectively, one needs to classify which dynamical
regime $\nu_j$ the system is in so that the correct ${\bf \Phi}^{\nu_j}$ can be used for the
reconstruction.  The classification is predicated on the idea that our measurement locations (\ref{eq:pi})
have the sparsest reconstruction (of the form (\ref{eq:reconstr})), given the correct
parameter $\nu_j$~\cite{bright13,brunton13}.  To employ this fact in a classification algorithm,
we shall represent the point ${\bf p}_{i}$ using all the possible parameters $\mu_j$,
generate an underdetermined system, and use a compressive sensing algorithm to obtain
the sparsest representation. Our hypothesis is that the non-zero coefficients
$A_{i_{k},j_{k}}^{\nu}$ and $B_{i_{k},j_{k}}^{\nu}$ will be concentrated
in the correct $\nu$. Specifically we shall construct an underdetermined
linear system and solve it using a sparsity promoting algorithm in order to identify the correct $\mu_j$.
Thus classification is based upon $\ell_1$ optimization whereas reconstruction relies on
an $\ell_2$ projection unto the correctly identified basis.

To be more precise, our classification is based on the following minimization problem
\begin{equation}
\min\left\Vert \mathbf{a}\right\Vert _{1}\mbox{ s.t. }\mathbf{p}=\Psi\mathbf{a}
\end{equation}
where the unknown vector ${\bf a}$ determines the weighting of the library elements ${\bf \Psi}$ unto the
measurement ${\bf p}$.  The $\ell_1$ minimization is a sparsity promoting algorithm that attempts to give
a solution ${\bf a}$ that has as many zeros as possible.  The non-zero components are used as an indicator
for the dynamical regime and modes, $\nu_j$ and ${\bf \Phi}^{\nu_j}$ respectively.
Here, however, the classification exploits both time and space measurements, whereas previous work~\cite{bright13,brunton13}
only considered measurements at a single fixed time for classification. 

Classification can be accomplished by considering a relaxed form of the minimization by a LASSO algorithm~\cite{lasso} 
\begin{equation}
\min\left\Vert \mathbf{a}\right\Vert _{1}\mbox{ s.t. }\left\Vert \mathbf{p}-\Psi\mathbf{a}\right\Vert _{2}<\delta,\label{eq:Lasso}
\end{equation}
for some parameter $\delta$.  This minimization problem minimizes the error, in an $\ell_2$ sense, between the measurements
and the projection on to the library with the sparsest vector ${\bf a}$.  The sparse solution ${\bf a}$ can then be used
to determine the dynamical regime $\nu_j$ the system is in.  Once determined, reconstruction can be achieved
via (\ref{eq:reconstr}).  Specifically, the classification of the dynamical regime is done by summing the absolute
value of the coefficients of ${\bf a}$ that corresponds to each dynamical regime $\nu_j$.
To account for regimes with a larger number of coefficients allocated for their dynamical
regime, we normalize by dividing by the square root of the number
of POD modes allocated in ${\bf a}$ for each $\mu_j$. The classified regime $\mu_j$
is the one that has the largest magnitude after execution of this algorithm.

\section{Application:  Flow Around a Cylinder}

We demonstrate the sparse time-space sampling algorithm developed on 
a canonical problem of applied mathematics:  the
flow around a cylinder.  This problem is well understood and has already been the subject of
studies concerning sparse spatial measurements~\cite{bright13,gappy,lionel,lionel2}.  Specifically, it is 
know that for low to moderate Reynolds numbers, the dynamics is spatially
low-dimensional and POD approaches have been successful in quantifying the
dynamics~\cite{gappy,pod1,pod2,pod3,pod4}.  Additionally, the time dynamics is nearly periodic, thus suggesting
that the asynchronous sampling in time advocated here can be exploited for this problem.
The Reynolds number plays the role of the bifurcation parameter $\nu$ in (\ref{eq:Nu}).

\begin{figure}[t]
\begin{centering}
\hspace*{-0.1in} \begin{overpic}[scale=0.6]{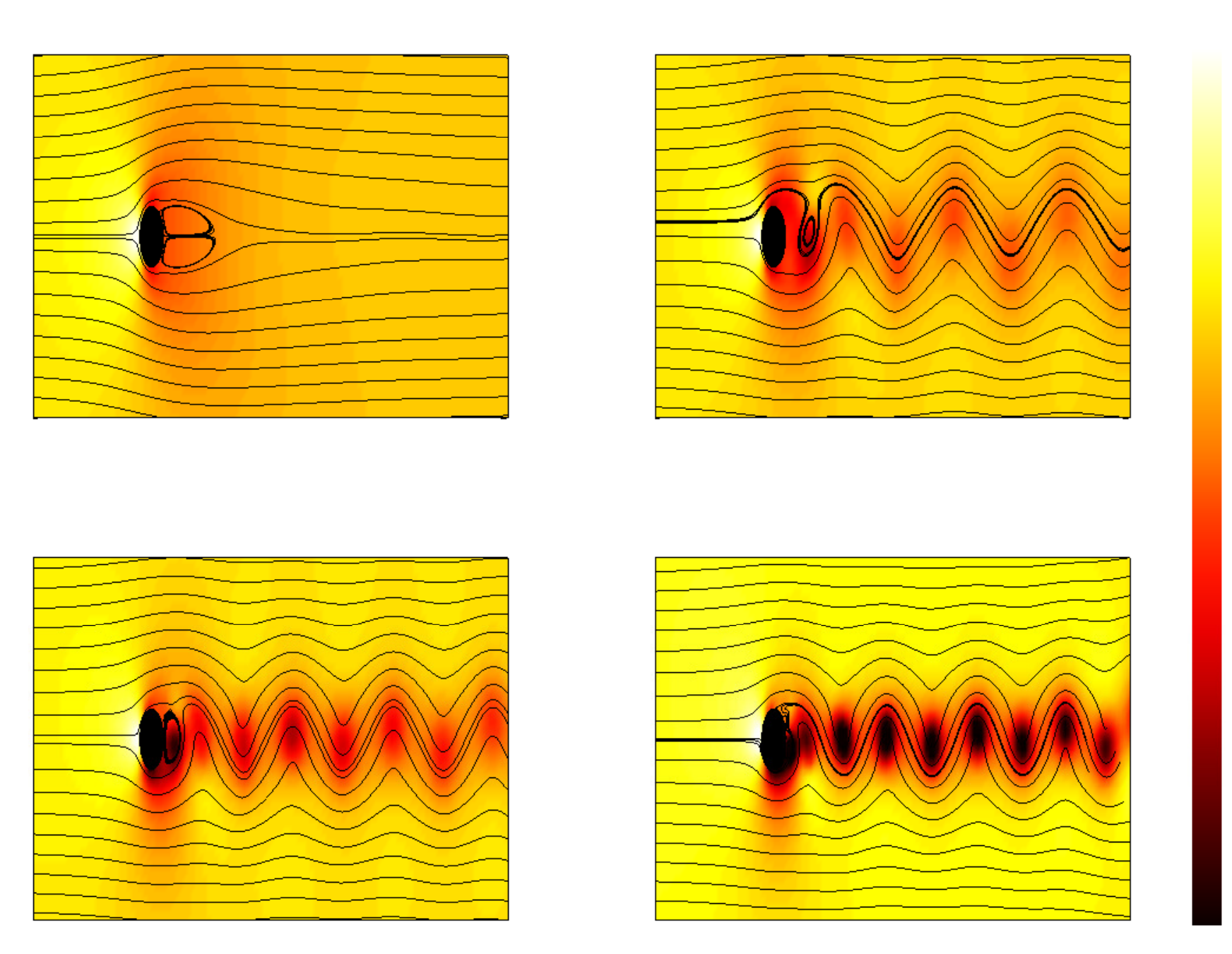}
\put(90,78){Pressure}
\put(102,72){0.5}
\put(102,25.3){-0.5}
\put(102,48.6){0.0}
\put(102,2){-1.0}
\put(19,75){$\nu=40$}
\put(18,34){$\nu=300$}
\put(69,75){$\nu=150$}
\put(67,34){$\nu=1000$}
\put(-1,2){-3}
\put(0,31){3}
\put(0,17){0}
\put(1,-1){-0.5}
\put(10,-1){0}
\put(25,-1){4}
\put(40,-1){8}
\put(32,-2){x}
\put(-3,24){y}
\end{overpic} 
\par\end{centering}
\caption{\label{fig:Stream-Lines} Stream lines and pressure field near the cylinder for a typical
snapshot taken for Reynolds numbers  (a) $\nu=40$, (b) $\nu=150$, (c) $\nu= 300$, and (d) $\nu=1000$.}
\end{figure}
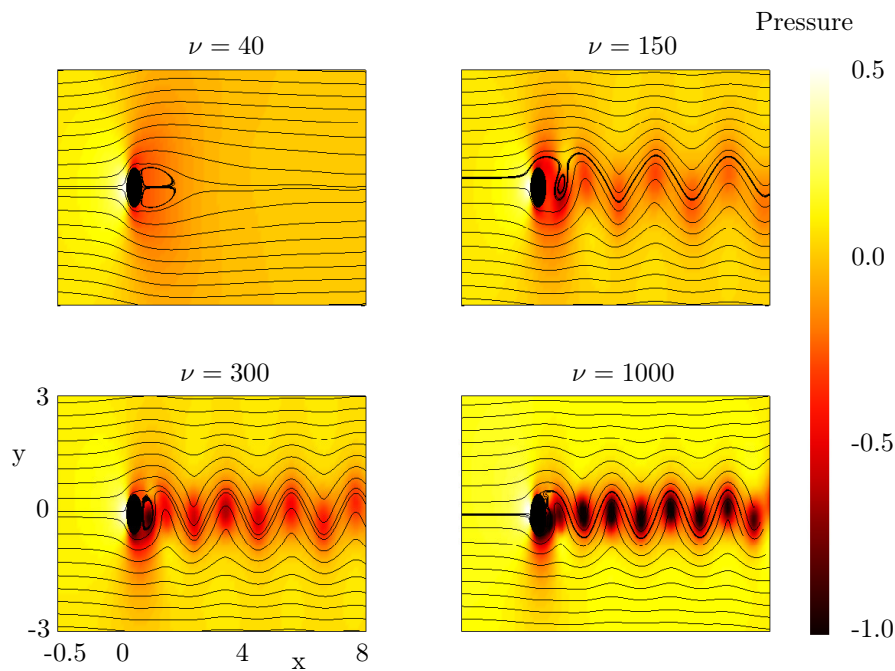

The data we consider comes from numerical simulations of the incompressible
Navier-Stokes equation: 
\begin{subeqnarray}
&&\frac{\partial u}{\partial t}+u\cdot\nabla u+\nabla p-\nu\nabla^{2}u=0\\
&&\nabla\cdot u=0
\label{eq:incompresNS}
\end{subeqnarray}
where $u\left(x,y,t\right)\in R^{2}$ represents the 2D velocity, and
$p\left(x,y,t\right)\in R^2$ the corresponding pressure field. The boundary condition are
as follows: (i) Constant flow of $u=\left(1,0\right)^{T}$ at $x=-15$, i.e., the
entry of the channel, (ii) Constant pressure of $p=0$ at $x=25$, i.e., the end of the channel, and (iii) 
Neumann boundary conditions, i.e. $\frac{\partial u}{\partial\mathbf{n}}=0$
on the boundary of the channel and the cylinder (centered at $(x,y)=(0,0$ and of radius unity).

\subsection{The Data:  Snapshots of Fluid Pressure}

Our data comes from numerical simulation of the incompressible Navier-Stokes
equation (\ref{eq:incompresNS}), computed using the Nektar++ package~\cite{nektar}, which produces accurate
results using a high-order finite element method. The algorithm uses a
non-uniform mesh consisting of 228 grid-points, with the high-order
method producing 66,000 Gaussian quadrature points, i.e. the numerical discretization yields an $n=66,000$ dimensional system.
Figure \ref{fig:Stream-Lines} shows
typical snapshots of the streamlines and pressure profiles taken
once the system transients have died away and the flow has reached 
equilibrium or a time-periodic state for various Reynolds numbers.  

The different dynamical regimes of (\ref{eq:incompresNS}) correspond to different
Reynolds number flow.  Thus the classification and reconstruction problem we consider 
involves library elements extracted from the flow at
Reynolds numbers $\nu=40,150,300,800,1000$.  Specifically, our measurements are
based on the pressure field at the perimeter of the cylinder where sensors
could be easily placed.  In particular, a sparse number of such sensors
are used to classify (identify) the Reynolds number and reconstruct the
pressure field around the cylinder and in the fluid flow.
Note that placing sensors in the fluid flow itself may be valid for
lab experiments and computer simulations.  However no realistic system
of interest, such as an airfoil or insect wing, places sensors directly
in the flow field behind the structure.  Thus there is a clear reliance on 
limited measurements on the cylinder itself.



\begin{figure}[t]
\begin{centering}
\vspace*{.2in}
\hspace*{-0.1in} \begin{overpic}[scale=0.6]{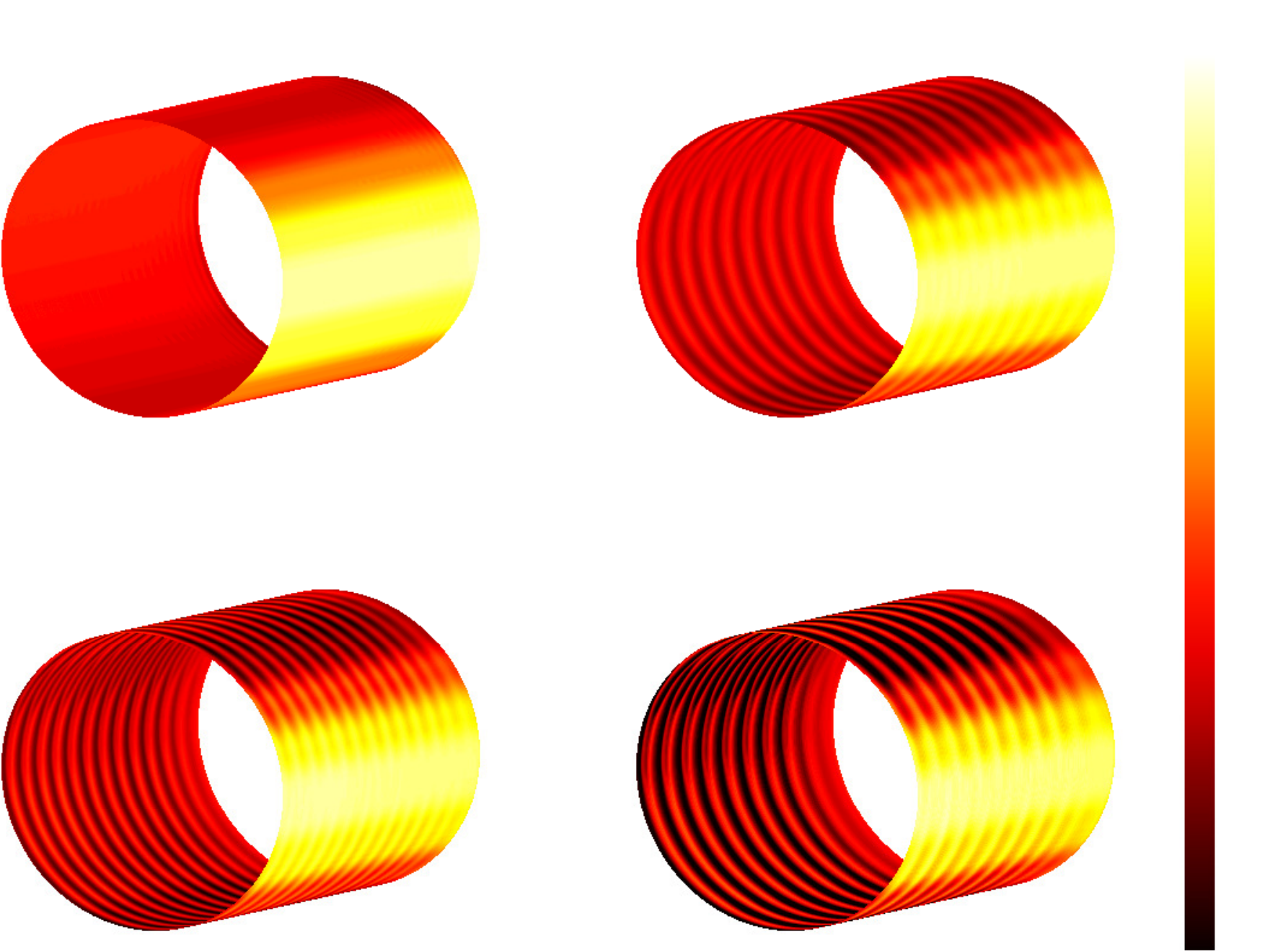}
\put(90,78){Pressure}
\put(102,72){0.5}
\put(102,25.3){-0.5}
\put(102,48.6){0.0}
\put(102,2){-1.0}
\put(29,71){$\nu=40$}
\put(28,31){$\nu=300$}
\put(79,71){$\nu=150$}
\put(77,31){$\nu=1000$}
\put(-3,14){0}
\put(9,28){$\pi/2$}
\put(25,14){$\pi$}
\put(9,-2){$3\pi/2$}
\put(34,40){time}
\put(15,35){0}
\put(28,37.5){$T$}
\put(15,38){\rotatebox{12}{$\xrightarrow{\hspace*{1.5cm}}$}}
\end{overpic} 
\par\end{centering}
\caption{\label{fig:3d plots}Time evolution of the pressure field around the cylinder for increasing
Reynolds numbers  $\nu=40$ ($T=8$), $\nu=150$ ($T=8$), $\nu= 300$ ($T=20$), and $\nu=1000$ ($T=20$).}
\end{figure}

The time evolution of the pressure field on the
perimeter of the cylinder is shown in Figure~\ref{fig:3d plots}
for several Reynolds numbers.  As can be seen, increasing the Reynolds
number drives the dynamics from a steady-state configuration to a
time-periodic evolution. The waterfall plots of the pressure dynamics
on the cylinder are essentially the snapshots used for projecting
the dynamics onto a low-dimensional manifold through the SVD, i.e. the 
snapshots provide the data necessary for the low-rank, POD methodology~\cite{PODHolmes}.

\begin{figure}[t]
\begin{center}
\vspace*{0.3in}
\begin{overpic}[scale=0.5]{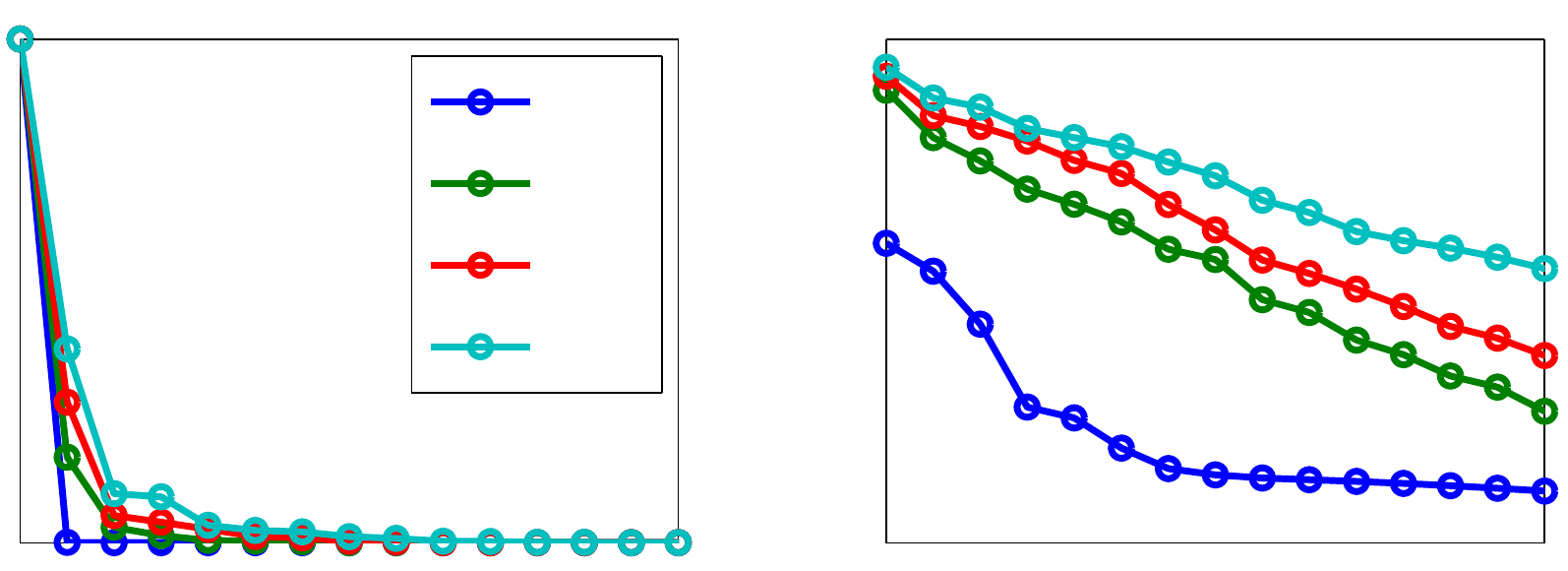}\end{overpic}
\begin{overpic}[scale=0.6]{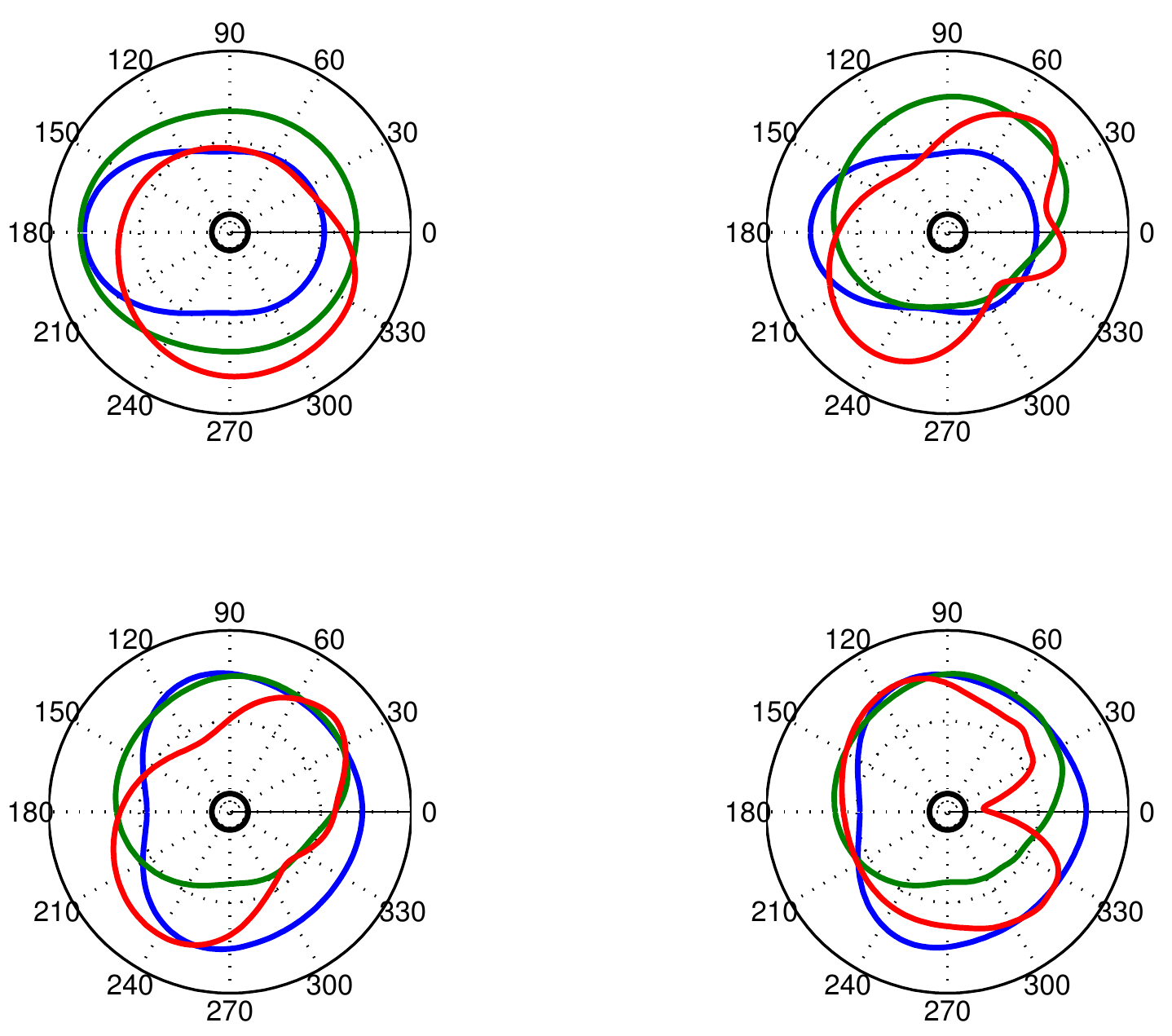}
\put(35,81){$\underline{{\nu=40}}$}
\put(34,33){$\underline{\nu=300}$}
\put(96,81){$\underline{\nu=150}$}
\put(95,33){$\underline{\nu=1000}$}
\put(20,116){$\nu\!=$}
\put(36,116){$40$}
\put(35,111){$150$}
\put(35,106){$300$}
\put(35,101){$1000$}
\put(53,120){0}
\put(52,110){-5}
\put(51,100){-10}
\put(51,89.5){-15}
\put(-5,88){\rotatebox{90}{Singular values ($\sigma_j^\nu/\sigma_1^\nu$)}}
\put(46,89){\rotatebox{90}{ $\log \left (\sum_{1}^{j} \sigma_k^\nu/ \sum_{1}^{N} \sigma_k^\nu \right)$}}
\put(-0.5,89.5){0}
\put(-0.5,120){1}
\put(2,87){1}
\put(13.5,87){5}
\put(27,87){10}
\put(42,87){15}
\put(56,87){1}
\put(67.5,87){5}
\put(81,87){10}
\put(96,87){15}
\linethickness{1.0mm}
\put(23,19){\textcolor{magenta}{{\line(1,0){12}}}}
\put(19.7,15.3){\textcolor{magenta}{-0.2}}
\put(33.1,15.3){\textcolor{magenta}{0.2}}
\put(38.5,14.5){\textcolor{magenta}{\bf{Pressure}}}
\put(23,19){\textcolor{magenta}{\circle*{2}}}
\put(35,19){\textcolor{magenta}{\circle*{2}}}
\put(19.5,19){\circle*{3}}
\put(81,19){\circle*{3}}
\put(19.5,68.5){\circle*{3}}
\put(81,68.5){\circle*{3}}
\put(40,52){\textcolor{blue}{\line(1,0){10}}}
\put(40,48){\textcolor{green}{\line(1,0){10}}}
\put(40,44){\textcolor{red}{\line(1,0){10}}}
\put(45,40){\circle*{3}}
\put(39.5,56){Dominant Modes}
\put(53,51.5){mode 1}
\put(53,47.5){mode 2}
\put(53,43.5){mode 3}
\put(53,39.5){cylinder}
\end{overpic} 
\end{center}
\vspace*{-.2in}
\caption{\label{fig:norm sing val} Normalized decay of singular values (top left) of the data matrix for various Reynolds number $\nu\!=\!40, 150, 300, 1000$  (top right on logarithmic scale). As the Reynolds number increases, a higher dimensional space is needed to preserve the same amount of energy.  The dominant three pressure modes are shown in
polar coordinates.  The pressure scale is in magenta (bottom left). }
\end{figure}
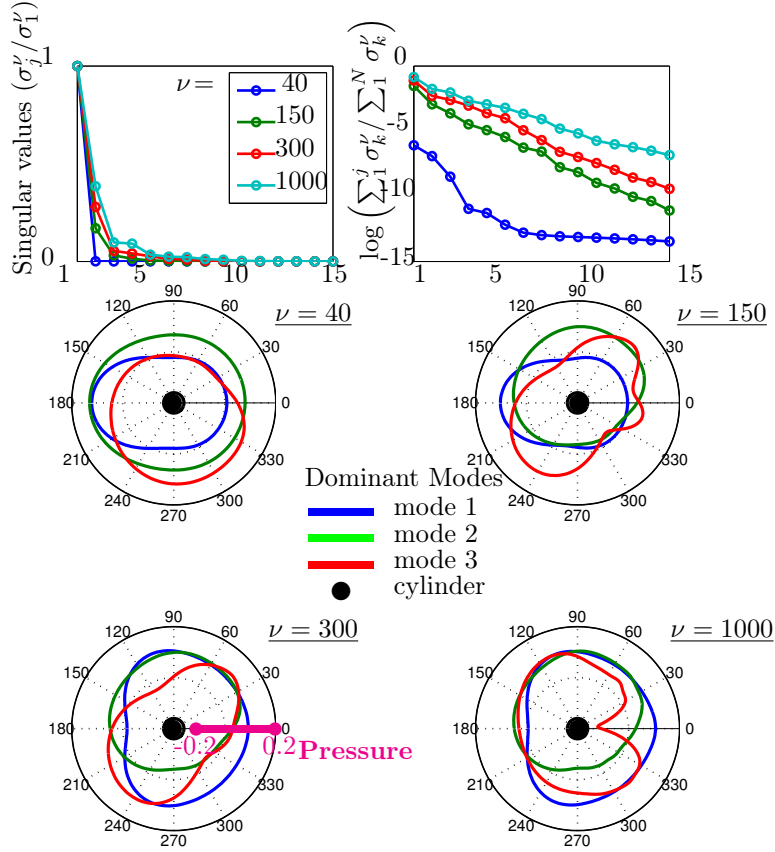

For each relevant value of the parameter $\nu$ we perform an SVD
on the data matrix in order to apply the dimensionality reduction method described
in Sec.~\ref{sec:Dimension-Reduction-Method}. It is well known
that for relatively low Reynolds number, a fast decay of the singular
values is observed. This can be seen in Fig.~\ref{fig:norm sing val} along with
the 3 most dominant POD modes.
Thus for every Reynolds number $\nu=40,150,300,800,1000$, we apply
the procedures described in Sec.~2 and 3 on the data matrix
${\bf A}^{\nu}$ in order to compute the quantities $\mathbf{v}_{i_{k}}^{\nu}$ and  $\mathbf{\omega}_{i_{k},j_{k}}^{\nu}$.


\subsection{Sparse Representation}

For every Reynolds number we keep 99\% of the total energy (variance) which gives for
the Reynolds numbers $\nu=40,150,300,800,1000$ considered a total of $1,3,3,9,9$ POD
modes to represent the dynamics.  More precisely, it determines the number of non-zero coefficients (1,3,3,9,9) 
in the sparse approximation \eqref{eq:sparse_app}
(i.e. non zero $\alpha_{k}^{c,\nu}$'s and $\alpha_{k}^{s,\nu}$'s).
The specific modes chosen by our representation, based on the singular
value analysis, namely on the magnitude of $\left(\alpha_{i,j}^{\nu}\right)^{2}=\left(\alpha_{i,j}^{c,\nu}\right)^{2}+\left(\alpha_{i,j}^{s,\nu}\right)^{2}$
is found in Table~\ref{ta:tab1}.

\begin{table}
\caption{Sparse Representation: The table contains the top 3 coefficients in the expansion \eqref{eq:sparse_app},
where we specify $\left(\alpha_{i,j}^{\nu}\right)^{2}=\left(\alpha_{i,j}^{c,\nu}\right)^{2}+\left(\alpha_{i,j}^{s,\nu}\right)^{2}$}
\begin{tabular}{|c||c|c|c|c|}
\hline 
Reynolds\\ Number & \# Modes & \multicolumn{3}{c|}{Top 3 Coefficients}\tabularnewline
\hline 
\hline 
$\nu=40$ & 1 & $\alpha_{1,0}^{40}=1$ &  & \tabularnewline
\hline 
\hline 
$\nu=150$ & 3 & $\alpha_{1,0}^{150}=0.97$ & $\alpha_{2,T=170}^{150}=0.03$ & \tabularnewline
\hline 
$\nu=300$ & 3 & $\alpha_{1,0}^{300}=0.92$ & $\alpha_{2,T=375}^{300}=0.07$ & \tabularnewline
\hline 
$\nu=800$ & 9 & $\alpha_{1,0}^{800}=0.87$ & $\alpha_{2,T=343}^{800}=0.11$ & $\alpha_{3,T=172}^{800}=0.04$\tabularnewline
\hline 
$\nu=1000$ & 9 & $\alpha_{1,0}^{1000}=0.85$ & $\alpha_{2,T=333}^{1000}=0.12$ & $\alpha_{3,T=166}^{1000}=0.06$\tabularnewline
\hline 
\end{tabular}
\label{ta:tab1}
\end{table}

\begin{figure}
\begin{center}
\begin{overpic}[scale=0.6]{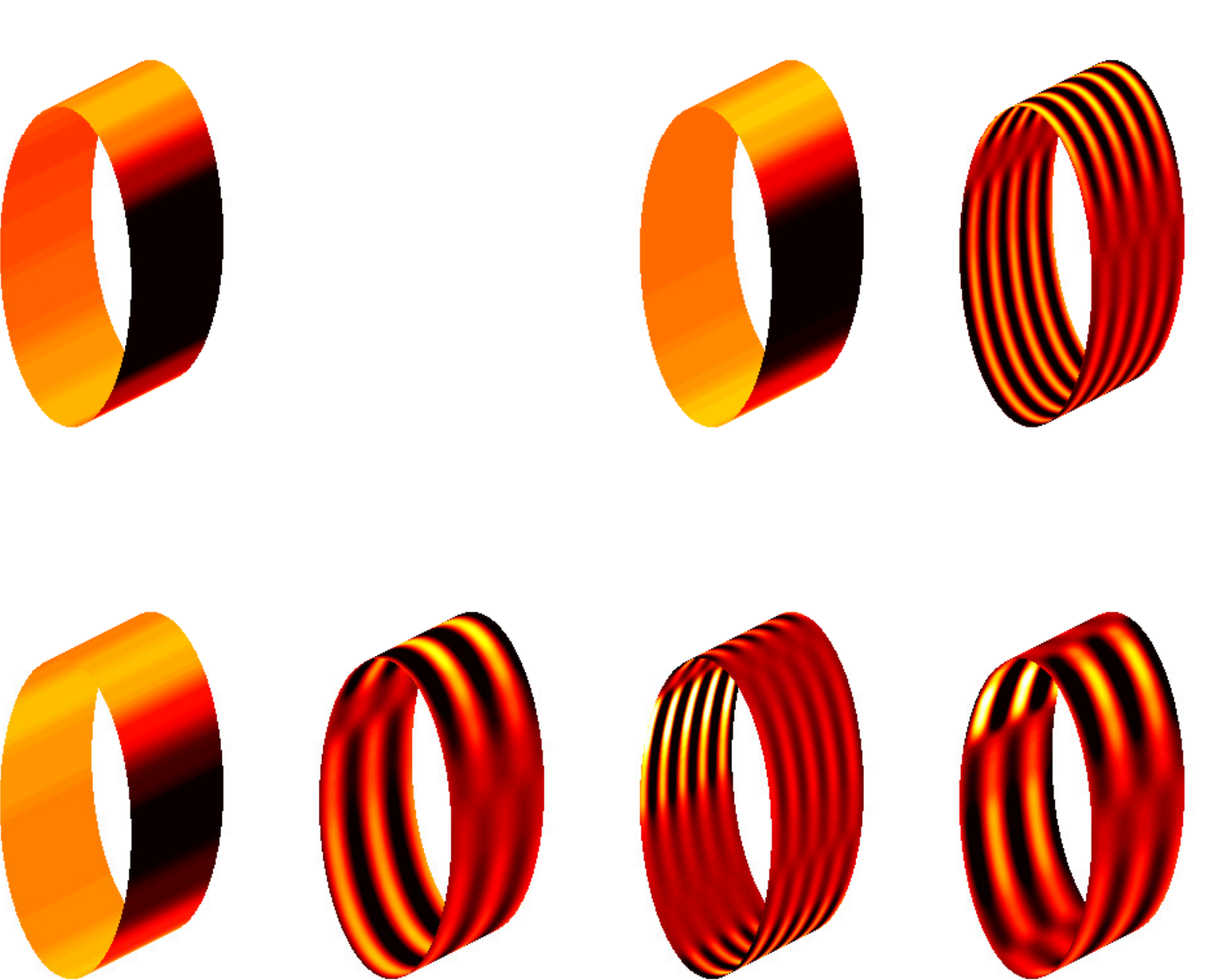}
\put(0,34){\underline{$\nu=800$}}
\put(0,79){\underline{$\nu=40$}}
\put(52,79){\underline{$\nu=150$}}
\put(-4.5,10){\rotatebox{90}{mode 1}}
\put(22,10){\rotatebox{90}{mode 2}}
\put(48.5,10){\rotatebox{90}{mode 3}}
\put(75,10){\rotatebox{90}{mode 4}}
\put(-4.5,55){\rotatebox{90}{mode 1}}
\put(48.5,55){\rotatebox{90}{mode 1}}
\put(75,55){\rotatebox{90}{mode 2}}
\put(5,39){0}
\put(15,43){$T$}
\put(5,42){\rotatebox{25}{$\xrightarrow{\hspace*{1.0cm}}$}}
\put(2,13){0}
\put(2,28){$\pi/2$}
\put(13,13){\textcolor{white}{$\pi$}}
\put(3,-3){$3\pi/2$}
\end{overpic}
\end{center}
\caption{\label{fig:Time Space Mods}The figure above shows the temporal-spatial
modes. These modes corresponds to the specified Reynolds number, and
the mode number corresponds to the magnitude of the energy of the
specific mode.}
\end{figure}
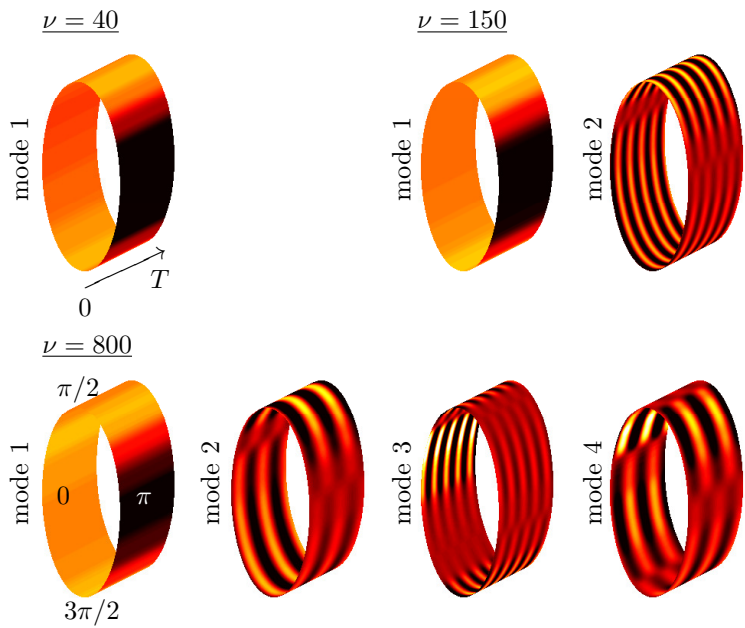

\subsection{Results}

We applied the sparse classification algorithm in Sec.~4 to
the training data collected in the previous subsection. Based on previous results \cite{bright13}, we
have chosen the sparse measurements at positions corresponding to
maxima and minima points of the dominant spatial modes (the
$\mathbf{v}_{i}^{\nu}$). Specifically 20 such position where chosen.
The classification algorithm was applied to 100 measurements randomly
selected over the time interval of training.  Namely, classification was based on 2.5\% of the data. The
results of classification subject to noise with a random distribution
is presented in Table~\ref{tab:Results:}.

\begin{table}
\caption{\label{tab:Results:}The error rate of our algorithm with noise
that is normally distributed with the standard deviation increasing across the columns.
It can be seen that for low to no additive error choosing a small
parameter $\delta$ for the Lasso \ref{eq:Lasso} provides a good performance,
while once the error is increased then the algorithm fails, and provides
result that are not any better than guessing. Note that increasing the Lasso parameter
$\delta$ provides improved performance for higher error, while decreasing
the performance for low or no error.}
\begin{centering}
\begin{tabular}{|c|c|c|c|c|}
\hline 
Error (std)  & $0$ & $2^{-7}$ & $2^{-5}$ & $2^{-4}$\tabularnewline
\hline 
\hline 
 $\delta=0.1$ & 0.89 & 0.88 & 0.2 & 0.2\tabularnewline
\hline 
$\delta=0.7$ & 0.8 & 0.82 & 0.83 & 0.82\tabularnewline
\hline 
\end{tabular}
\par\end{centering}
\end{table}

It should be noted once again that the method demonstrated here extends previous
work~\cite{bright13,brunton13} by considering the full spatio-temporal dynamics.  Specifically,
previous findings were associated with successfully classifying and reconstructing a single snapshot in
time.  Here, the much more difficult task is considered of reconstruction the entire spatio-temporal
behavior with limited measurements.  The results presented show the method is quite effective in this
task, thus providing a potentially transformational paradigm for taking advantage of limited
measurements of a complex, spatio-temporal system.

\section{Conclusion and Outlook}

Modern data methods and analysis techniques look to exploit underlying low-dimensional
structures in complex data sets.  In the case of time-dependent data generated from
complex systems, effective methods must capitalize simultaneously on the underlying
low-rank structures in both time and space.  In doing so, we have demonstrated that
one can use a sparse number of spatial sensors that sample randomly (asynchronously) in time to accurately classify
the dynamical state of the system even in the presence of noise.  The method advocated integrates
two transformative tools of analysis:  (i) model reduction (machine learning) methods for complex
systems and (ii) compressive sensing/sparse representation 

The current mathematical architecture and innovation extends previous work on sparse sampling
of complex systems by exploiting the nearly periodic time behavior known to be exhibited by
many complex systems, including the flow around a cylinder example considered here.  Indeed,
taking advantage of the sparse signal representation in time allows for better classification of
the dynamical state.  Moreover it equally, and maximally, exploits the low-dimensional behavior 
in both time and space.  Previous analysis relied on a single measurement in time for performing
such classification tasks, thus only exploiting the low-dimensional spatial structures.

The combination of dimensionality reduction and sparsity promoting techniques advocated here
can be applied in an exceptionally broad context. Indeed, such algorithmic
strategies can be used to enhance computation and efficiently identify
measurement locations in a given system by {\em remembering} the
system's key characteristics in both time and space. The low-rank spatio-temporal nature of the approximations
used allows construction of control algorithms that wrap around the
dimensionality reduction and sparsity infrastructure\textemdash{}another
appealing aspect of the methodology.  In conclusion, we advocate a general theoretical framework for 
complex systems whereby low-rank, spatio-temporal libraries representing
the optimal modal basis in time and space are constructed, or learned, from snapshot sampling of the dynamics.

\bibliographystyle{unsrt}
\bibliography{cit}

\end{document}